\documentclass{article}

\usepackage{PRIMEarxiv}

\usepackage[utf8]{inputenc} 
\usepackage[T1]{fontenc}    
\usepackage{hyperref}       
\usepackage{url}            
\usepackage{booktabs}       
\usepackage{amsfonts}       
\usepackage{nicefrac}       
\usepackage{microtype}      
\usepackage{lipsum}
\usepackage{fancyhdr}       
\usepackage{graphicx}       
\graphicspath{{media/}}     
\usepackage{fullpage}
\usepackage{graphicx}
\usepackage{multicol}
\usepackage[utf8]{inputenc}
\usepackage[english]{babel}
\usepackage{amsmath}
\usepackage{tikz}
\usepackage{pgfplots}
\usepackage{amssymb, graphics, amsmath}
\usepackage{multirow}
\usepackage{array}
\usepackage{bigints}
\usepackage{amsthm}
\usepackage{comment}
\usepackage{verbatim}
\usepackage{amsmath}

\newtheorem{theorem}{Theorem}[section]
\newtheorem{lemma}{Lemma}[section]
\newtheorem{result}{Result}[section]
\numberwithin{equation}{section}

\newtheorem{remark}{Remark}[section]
\newtheorem*{theorem-non}{Theorem}

\title{On the Law of Large Numbers and Convergence Rates for the Discrete Fourier Transform of Random Fields
}

\author{
  Vishakha \\
  Department of Mathematical Sciences \\
  University of Cincinnati \\
  Cincinnati\\
  \texttt{sharmav4@mail.uc.edu} \\
}

\begin{document}
\maketitle

\begin{abstract}
We study the $\text{Marcinkiewicz-Zygmund}$ strong law of large numbers for the cubic partial sums of the discrete Fourier transform of random
fields. We establish Marcinkiewicz–Zygmund types rate of convergence for the discrete Fourier transform of random
fields under weaker conditions than identical distribution.
\end{abstract}

\keywords{ Discrete Fourier transform \and Law of large numbers \and
Rate of convergence}

\section{Introduction}
For $d\geq1$, let $\mathbb{Z}_+^d$ denote the set of positive integer $d-$dimensional lattice points. Let `$\preceq$' be a component-wise partial order on $\mathbb{Z}_+^d$. For $\textbf{m}=(m_1, \dots, m_d),~ \textbf{n}=(n_1, \dots, n_d) \in \mathbb{Z}_+^d$, $\textbf{m}\preceq \textbf{n}$ means that $m_i \leq n_i$, $\forall 1\leq i\leq d$; $|\textbf{n}|$ is used to denote the product $n_1 n_2 \cdots n_d$, and $\textbf{n} \rightarrow \infty$ means that  $\min\{n_1, \dots, n_d\} \rightarrow \infty$. 

Suppose $\{X_{\textbf{n}}:{\textbf{n}}\in \mathbb{Z}_+^d\}$ and $X$ are real-valued random variables on some probability space, where $X$ has the same distribution as $X_\textbf{1}$ and  $EX_{\textbf{n}}=0$.
Throughout this paper, the notation $\log^+x$ denotes $\max\{1, \log x\}$. \\

For $t=(t_1, \dots, t_d) \in \mathbb{T}^d$, where $\mathbb{T}$ denotes $[-\pi,\pi)$, the discrete Fourier transform for random fields is given as
$S_{\textbf{n}}(t)= \sum_{\textbf{k}\preceq \textbf{n}} e^{i\textbf{k}\cdot t}X_{\textbf{k}},$
where $\textbf{k}\cdot t=k_1t_1+k_2t_2+\cdots +k_dt_d$.\\

It has been a keen interest to study the law of large numbers for the random variables and their convergence rates. A lot of significant research has been done in this direction by many authors for both $d=1$ and $d\geq 2$.
For the case $d=1$, both strong and weak laws hold under the assumption of a finite mean, but for $d\geq 2$, only the weak law holds under this assumption. Therefore, it was of interest to study strong laws for the higher dimension. \\

Many results have been obtained for the sequence of independent and identically distributed random variables.

For $d=1$, Marcinkiewicz \textit{et.al.} in \cite{Marcinkiewicz1937} proved the classical Marcinkiewicz–Zygmund strong
law of large numbers, that is, $n^{-1/p}S_n \xrightarrow{\text{a.s.}} 0$ if and only if  $E|X|^p< \infty$ with $EX=0$ and $1 \leq p <2$. 

For $d\geq 1$, Smythe \cite{smythe1974sums} showed that the strong law of large numbers holds if and only if $E|X|(\log^+|X|)^{d-1}$ is finite.  
Later Gut \cite{gut1978marcinkiewicz} generalized this result and proved the strong Marcinkiewicz law, that is, $|\textbf{n}|^{-1/p}S_{\textbf{n}} \xrightarrow{\text{a.s.}} 0$ if and only if  $E|X|^p(\log^+|X|)^{d-1}< \infty$ with $EX=0$ and $1 \leq p <2$. \\

It has been fascinating to investigate the convergence rates in the law of large numbers, that is, to investigate how fast quantities such as $P\big(|S_{\textbf{n}}|>\epsilon|\textbf{n}|^{\alpha}\big)$ approach to zero under some sort of finite moment conditions. In order to study these quantities, one can study the convergence of sums of type   $\sum_{\textbf{n}=1}^{\infty}|\textbf{n}|^{\beta}P\big(|S_{\textbf{n}}|>\epsilon|\textbf{n}|^{\alpha}\big)$, where $\alpha$ and $\beta$ are real numbers.
Several authors such as Baum \cite{baum1965convergence}, Katz \cite{katz1963probability}, Marcinkiewicz \cite{Marcinkiewicz1937} \cite{gut1978marcinkiewicz} etc. have done a significant work in investigating the convergence of these sums for a  sequence of i.i.d. random variables. 

In the case $d=1$ and $\alpha=0$, the initial works were carried out by Hsu and Robbins \cite{hsu1947complete} and Erdos \cite{erdos1949theorem} for $\beta=0$, and by Spitzer \cite{spitzer1956combinatorial} for $\beta=-1$, and later further investigations were continued by many other researchers.\\

For higher dimensions, $d\geq 2$, 
Symthe \cite{smythe1974sums} and others studied the sum $\sum_{\textbf{n}=1}^{\infty}|\textbf{n}|^{\beta}P\big(|S_{\textbf{n}}|>\epsilon|\textbf{n}|\big)$ for $\beta=0$ and $\beta \geq -1$, respectively.\\

For a sequence of pairwise independent and identically distributed random variables several authors, viz., Etemadi \cite{etemadi1981elementary}, Baum and Katz \cite{baum1965convergence} etc. have contributed significantly in studying the law of large numbers and their rates of convergence in the case $d=1$. 
For a sequence of random
variables that are not necessarily independent and identically distributed, Zhang in \cite{zhang2017law} proved the strong law of large numbers for the discrete Fourier transform of $(X_n)_n$ with finite $pth$ moment for $p=1$ and their Marcinkiewicz–Zygmund type rate of convergence for $1 < p < 2$ under
weaker conditions than identical distribution. In the work of \cite{zhang2017law}, the convergence rates were left unstudied for the case $p=1$. \\

In the same direction, our goal in this paper is to generalize certain results of \cite{zhang2017law} to the discrete Fourier transform of random fields. Under weaker conditions than identical distribution, we specifically study the convergence of $\frac{S_{\textbf{n}}(t)}{|\textbf{n}|}$ in the $L_1$-norm, and the convergence of sums of type 
$\sum_{\textbf{n}=\textbf{1}}^{\infty}|\textbf{n}|^{p/r-2} P\big(|S_{\textbf{n}}(t)|>\epsilon|\textbf{n}|^{1/r}\big)$, 
and for the restricted index set, sums of type $\sum_{\textbf{n}=\textbf{1}}^{\infty}|\textbf{n}|^{p/r-1-1/d} P\big(\underset{\textbf{1}\preceq \textbf{m} \preceq \textbf{n}}{\max}
|S_{\textbf{m}}(t)|>\epsilon|\textbf{n}|^{1/r}\big)$, where $1\leq p<2$.\\

The paper is organized as follows: in Section \ref{stateMainResults}, we state our main results. In Section \ref{lemmaWithProofs}, we collect certain lemmas and their proofs. Finally, the main results are proved in Section \ref{theoemproofs}.\\

\section{Notations}
The following notations are used throughout the paper.

\begin{itemize}
     \item Define  $d(x)=Card\{\textbf{n}\in \mathbb{Z}_+^d : |\textbf{n}| = [x]\}$ and $M(x)=\sum_{j=1}^{[x]}d(j)$.
    \item For any $A,B \geq 0$, $A\ll B$ means that there exists a constant $c>0$ such that $A\leq cB$.
    \item For $a_n,b_n>0$, $\sum_{n=0}^{\infty}a_n \sim \sum_{n=0}^{\infty}b_n$ means $\frac{a_n}{b_n}\rightarrow 1$.
    \item Symbols $\lambda$ and $L$ denote the lebesgue measure on $\mathbb{T}$ and the product measure, $\lambda^d$, on $\mathbb{T}^d$, respectively. 
    \item $\tilde{\Omega}:=\Omega\times \mathbb{T}^d$ denotes the product space and $\tilde{P}:=P\times L$ denotes the product measure on this space.
\end{itemize}

\section{Main results}\label{stateMainResults}

For $\textbf{1} \preceq \textbf{n}$, let $X_{\textbf{n}}$ and $X$ be real-valued random variables on a probability space $(\Omega,\mathcal{F},P)$. 
No dependence between variables is assumed.\\

We state our main results in this section and they will be proved later in Section \ref{theoemproofs}. \\

\begin{theorem} \label{Theorem1}
Assume $(X_{\textbf{n}})_{\textbf{n}}$ satisfies condition
    \begin{equation}
    \frac{1}{|\textbf{n}|}\sum_{\textbf{k}\preceq \textbf{n}}P(|X_{\textbf{k}}|> x) \leq MP(|X|> x), \text{ for all } x\geq 0,~ 1\preceq \textbf{n}. \tag{MD}
\end{equation}
If $E|X| < \infty$, then 
$$\lim_{\max\{n_1,\dots,n_d\}\rightarrow \infty}\frac{S_{\textbf{n}}(t)}{|\textbf{n}|}=0 \quad \text{in the }L_1(\tilde{\Omega})\text{-norm}.$$\\
\end{theorem}

\begin{theorem} \label{Theorem2}
Assume $(X_{\textbf{n}})_{\textbf{n}}$ satisfies condition (MD).
Let $1<p<2$, $1\leq r\leq p$. \\

If $E|X|^p(\log^+|X|)^{d-1} < \infty$, then for every $\epsilon >0$ and for almost all $t\in \mathbb{T}^d$
$$\sum_{\textbf{n}=1}^{\infty}|\textbf{n}|^{p/r-2} P\big(|S_{\textbf{n}}(t)|>\epsilon|\textbf{n}|^{1/r}\big) < \infty.$$\\   
\end{theorem}

\noindent For the case $p=1$, the above Theorem \ref{Theorem2} holds under following additional moment condition.\\

\begin{theorem}\label{Theorem3}
Assume $(X_{\textbf{n}})_{\textbf{n}}$ satisfies condition (MD).\\

If $E\big[|X|(\log^+|X|)^{d}\big] < \infty$, then for every $\epsilon >0$ and for almost all $t\in \mathbb{T}^d$
$$\sum_{\textbf{n}=1}^{\infty}|\textbf{n}|^{-1} P\big(|S_{\textbf{n}}(t)|>\epsilon|\textbf{n}|\big) < \infty.$$\\    
\end{theorem}

\noindent
Denote $M=(m,\dots,m)$, $N=(n,\dots,n)\in \mathbb{Z}_{+}^d.$\\

\begin{theorem}\label{Theorem4}
Assume $(X_{\textbf{n}})_{\textbf{n}}$ satisfies condition (MD).
Let $1<p<2$ and $1\leq r\leq p$.\\

If $E|X|^p < \infty$, then for every $\epsilon >0$ and for almost all $t\in \mathbb{T}^d$
$$\sum_{n=1}^{\infty}|N|^{p/r-1-1/d} P\bigg(\underset{1\leq m \leq n}{\max}
|S_{M}(t)|>\epsilon|N|^{1/r}\bigg) < \infty.$$    
\end{theorem}


\begin{remark}\label{remark1}
The assumptions of Theorem \ref{Theorem4} imply the following Marcinkiewicz–Zygmund strong law of large numbers for the cubic sums.\\

For almost all $t\in \mathbb{T}^d$
$$\lim_{N \rightarrow \infty}\frac{S_{N}(t)}{|N|^{1/p}}=0, \quad \text{P-a.s.}$$  \\ 
\end{remark}

\noindent For the case $p=1$, the above Theorem \ref{Theorem4} holds under following additional moment condition.

\begin{theorem}\label{Theorem5}
Assume $(X_{\textbf{n}})_{\textbf{n}}$ satisfies condition (MD).\\

If $E|X|log^+|X| < \infty$, then for every $\epsilon >0$ and for almost all $t\in \mathbb{T}^d$
$$\sum_{n=1}^{\infty}|N|^{-1/d} P\bigg(\underset{1\leq m \leq n}{\max}
|S_M(t)|>\epsilon|N|\bigg) < \infty.$$\\   
\end{theorem}

\begin{remark}
The assumptions of Theorem \ref{Theorem5} imply that for almost all $t\in \mathbb{T}^d$
$$\lim_{n \rightarrow \infty}\frac{S_{N}(t)}{|N|}=0, \quad \text{P-a.s.}$$ \\  
\end{remark}

\begin{remark}
The analogues of Theorem \ref{Theorem4} and Theorem \ref{Theorem5} do not hold for the rectangular
partial sums since the maximal inequality of Theorem 4.4 \cite{weisz2012summability} does not always hold for rectangular
partial sums (see Fefferman \cite{fefferman1971divergence}).   
\end{remark}

\section{Lemmas}\label{lemmaWithProofs}

We need following lemmas in order to prove our main results.\\

\begin{lemma}\label{lemma2}
Assume $(X_{\textbf{n}})_{\textbf{n}}$ satisfies condition (MD). Let $1\leq p<2$, $1\leq r\leq p$.\\

If $E\big[|X|^p(\log^+|X|)^{d-1}\big] < \infty$, 
then 
$$(i) 
\sum_{\textbf{n}=1}^{\infty}|\textbf{n}|^{(p-1)/r-2}\sum_{\textbf{k} \preceq \textbf{n}}E\big[|X_{\textbf{k}}|I(|X_{\textbf{k}}|>|\textbf{n}|^{1/r})\big] < \infty, \quad \text{for } p>1; 
$$
and,
$$(ii)
\sum_{\textbf{n}=1}^{\infty}|\textbf{n}|^{(p-2)/r-2}\sum_{\textbf{k} \preceq \textbf{n}}E\big[|X_{\textbf{k}}|^2I(|X_{\textbf{k}}|\leq|\textbf{n}|^{1/r})\big] < \infty, \quad \text{for } p\ge1.
$$\\

\begin{proof}

\textbf{(i)} Using condition (MD) and a known fact (\ref{EXuppertail}), we obtain
\begin{equation*}
\begin{split}
\sum_{\textbf{k} \preceq \textbf{n}}E\big[|X_{\textbf{k}}|I(|X_{\textbf{k}}|>|\textbf{n}|^{1/r})\big]
&=|\textbf{n}|^{1/r}\sum_{\textbf{k} \preceq \textbf{n}}P(|X_{\textbf{k}}|>|\textbf{n}|^{1/r})
+\bigintsss_{|\textbf{n}|^{1/r}}^{\infty}\sum_{\textbf{k} \preceq \textbf{n}}P(|X_{\textbf{k}}|>x)\,dx\\
&\leq M|\textbf{n}|^{1/r+1}P(|X|>|\textbf{n}|^{1/r})+M|\textbf{n}|\bigintsss_{|\textbf{n}|^{1/r}}^{\infty}P(|X|>x)\,dx.\\
\end{split}   
\end{equation*}

Then,
\begin{align}
\sum_{\textbf{n}=1}^{\infty}|\textbf{n}|^{(p-1)/r-2}\sum_{\textbf{k} \preceq \textbf{n}}E\big[|X_{\textbf{k}}|I(|X_{\textbf{k}}|>|\textbf{n}|^{1/r})\big]
&\ll\sum_{\textbf{n}=1}^{\infty}\bigg[|\textbf{n}|^{p/r-1}P(|X|>|\textbf{n}|^{1/r})+|\textbf{n}|^{(p-1)/r-1}\bigintssss_{|\textbf{n}|^{1/r}}^{\infty}P(|X|>x)\,dx\bigg]\notag \\ 
&\leq \sum_{\textbf{n}=1}^{\infty}|\textbf{n}|^{p/r-1}P(|X|>|\textbf{n}|^{1/r})
+\sum_{\textbf{n}=1}^{\infty}|\textbf{n}|^{(p-1)/r-1}\bigintssss_{|\textbf{n}|^{1/r}}^{\infty}P(|X|>x)\,dx\notag \\ 
&=\sum_{j=1}^{\infty}j^{p/r-1}d(j)P(|X|>j^{1/r})
+\sum_{j=1}^{\infty}j^{(p-1)/r-1}d(j)\bigintssss_{j^{1/r}}^{\infty}P(|X|>x)\,dx \notag \\
&=(I)+(II).\label{L4.2E1}
\end{align}

As a consequence of Lemma 2.1 of \cite{gut1978marcinkiewicz}, we obtain that series
$(I)<\infty$ under the given moment condition.\\  

To evaluate $(II)$, we start by writing $\bigintsss_{j^{1/r}}^{\infty}P(|X|>x)\,dx=\sum_{i=j}^{\infty}b_i$, where $b_i=\bigintsss_{i^{1/r}}^{(i+1)^{1/r}}P(|X|>x)\,dx$. So, by changing the order of summation
$(II)
=\sum_{i=1}^{\infty}b_i\sum_{j=1}^{i}j^{(p-1)/r-1}d(j)$.\\

Now, since $d(j)=M(j)-M(j-1)$,  we have
\begin{align}
(II) \ll \sum_{i=1}^{\infty}b_i\sum_{j=1}^{i-1}M(j)~j^{(p-1)/r-2} ~ +~ \sum_{i=1}^{\infty}b_i~i^{(p-1)/r-1}M(i)
:=J_1+J_2. \label{L4.2E2}    
\end{align}
\\

Let us first consider $J_2$. Since we have the relation $M(i) \sim i(\log i)^{d-1}$ as $i \rightarrow \infty$, therefore
\begin{align}
J_2 \sim \sum_{i=1}^{\infty}b_i~i^{(p-1)/r}(\log i)^{d-1} 
\ll \sum_{i=1}^{\infty}\bigintssss_{i^{1/r}}^{(i+1)^{1/r}}x^{p-1}(\log^+ x)^{d-1}P(|X|>x)\,dx. \label{L4.2E3}
\end{align}
\\

Let us now consider $J_1$. By using the relation between $M(j)$ and $\log j$ as we used above, we get
\begin{align}
J_1\ll \sum_{j=1}^{\infty}M(j)j^{(p-1)/r-2}\sum_{i=j}^{\infty}b_i 
&\sim \sum_{j=1}^{\infty}j^{(p-1)/r-1}(\log j)^{d-1}\sum_{i=j}^{\infty}b_i. \label{L4.2E4}
\end{align}
By monotonicity, for $i\geq j$
\begin{align}
(\log j)^{d-1}b_i \leq \bigintssss_{i^{1/r}}^{(i+1)^{1/r}}(\log^+ x)^{d-1}P(|X|>x)\,dx. \label{L4.2E5}  
\end{align}

Then, by using (\ref{L4.2E5}) and doing the change of order of summation, 
\begin{align}
\sum_{j=1}^{\infty}j^{(p-1)/r-1}(\log j)^{d-1}\sum_{i=j}^{\infty}b_i
\ll \sum_{i=1}^{\infty} i^{(p-1)/r} \bigintssss_{i^{1/r}}^{(i+1)^{1/r}}(\log^+ x)^{d-1}P(|X|>x)\,dx \notag \\
\leq \sum_{i=1}^{\infty}\bigintssss_{i^{1/r}}^{(i+1)^{1/r}}x^{p-1} (\log^+ x)^{d-1}P(|X|>x)\,dx. \label{L4.2E6}
\end{align}
\\

By combining (\ref{L4.2E3}) together with (\ref{L4.2E4}) and (\ref{L4.2E6}), we obtain that $J_1+J_2$ is finite (and, so is $(II)$) by our given moment condition $E\big[|X|^p(\log^+|X|)^{d-1}\big] < \infty$.\\

Therefore, we conclude
$$\sum_{\textbf{n}=1}^{\infty}|\textbf{n}|^{(p-1)/r-2}\sum_{\textbf{k} \preceq \textbf{n}}E\big[|X_{\textbf{k}}|I(|X_{\textbf{k}}|>|\textbf{n}|^{1/r})\big] < \infty.$$\\

\noindent Next, we prove the second part of the lemma.\\

\textbf{(ii)} Using condition (MD) and a known fact (\ref{EXlowertail}), we obtain
\begin{align}
\sum_{\textbf{k} \preceq \textbf{n}}E\big[|X_{\textbf{k}}|^2I(|X_{\textbf{k}}|\leq|\textbf{n}|^{1/r})\big]
\leq \sum_{\textbf{k} \preceq \textbf{n}}\bigintssss_{0}^{|\textbf{n}|^{2/r}}P(X_{\textbf{k}}^2>x)\,dx
&\leq M|\textbf{n}|\bigintssss_{0}^{|\textbf{n}|^{2/r}}P(X^2>x)\,dx. \label{L4.2.2E1}
\end{align}
\\

So, 
\begin{align}
\sum_{\textbf{n}=1}^{\infty}|\textbf{n}|^{(p-2)/r-2}\sum_{\textbf{k} \preceq \textbf{n}}E\big[|X_{\textbf{k}}|^2I(|X_{\textbf{k}}|\leq|\textbf{n}|^{1/r})\big]
&\ll \sum_{\textbf{n}=1}^{\infty}|\textbf{n}|^{(p-2)/r-1}\bigintssss_{0}^{|\textbf{n}|^{2/r}}P(X^2>x)\,dx.  \label{L4.2.2E2}  
\end{align}
\\

Now, by writing $\bigintss_{0}^{|\textbf{n}|^{2/r}}P(X^2>x)\,dx 
= \sum_{i=1}^{j}\bigintss_{(i-1)^{2/r}}^{i^{2/r}}P(X^2>x)\,dx$, and doing the change of order of summation, (\ref{L4.2.2E2}) gives
\begin{align}
\sum_{\textbf{n}=1}^{\infty}|\textbf{n}|^{(p-2)/r-2}\sum_{\textbf{k} \preceq \textbf{n}}E\big[|X_{\textbf{k}}|^2I(|X_{\textbf{k}}|\leq|\textbf{n}|^{1/r})\big]
&\ll \sum_{i=1}^{\infty}\bigg[\sum_{j=i}^{\infty}j^{(p-2)/r-1} d(j)\bigg]\bigintssss_{(i-1)^{2/r}}^{i^{2/r}}P(X^2>x)\,dx. \label{L4.2.2E3}
\end{align}
\\

One can now easily verify that
\begin{align}
\text{ for any } \delta>0, \text{the convergence of } \sum_{j=i}^{\infty}j^{-\delta} d(j) \text{ is  equivalent to that of } \sum_{j=i}^{\infty}j^{-\delta-1} M(j). \label{L4.2.2E4}
\end{align}
and
\begin{align}
 \forall i\geq 3, \log i \leq \log (i-1)^2, \text{ and therefore}, \forall x\geq (i-1)^{2/r}, (\log i)^{d-1} \ll (\log^+x)^{d-1}.\label{L4.2.2E5}
\end{align} 
\\

By (\ref{L4.2.2E4}) and (\ref{L4.2.2E5}), we obtain
\begin{align}
\sum_{i=1}^{\infty}\bigg[\sum_{j=i}^{\infty}j^{(p-2)/r-1} d(j)\bigg]\bigintssss_{(i-1)^{2/r}}^{i^{2/r}}P(X^2>x)\,dx
& \sim \sum_{i=3}^{\infty}i^{\frac{(p-2)}{r}}(\log i)^{d-1}\bigintssss_{(i-1)^{2/r}}^{i^{2/r}}P(X^2>x)\,dx. \label{L4.2.2E6}
\end{align}
\\

Finally, by combining (\ref{L4.2.2E3}) together with (\ref{L4.2.2E6}), we  conclude that 
\begin{align*}
\sum_{\textbf{n}=1}^{\infty}|\textbf{n}|^{(p-2)/r-2}\sum_{\textbf{k} \preceq \textbf{n}}E\big[|X_{\textbf{k}}|^2I(|X_{\textbf{k}}|\leq|\textbf{n}|^{1/r})\big]  < \infty,
\end{align*}
by our moment condition
$E\big[|X|^p(\log^+|X|)^{d-1}\big] < \infty$.\\

This completes the proof of Lemma \ref{lemma2}.\\
\end{proof}
\end{lemma}

\vspace{0.5cm}

\noindent For the case $p=1$, the first part of Lemma \ref{lemma2} holds under the following additional moment condition.\\

\begin{lemma}\label{lemma3}
Assume $(X_{\textbf{n}})_{\textbf{n}}$ satisfies condition (MD)
and $E\big[|X|(\log^+|X|)^{d}\big] < \infty$. 
Then 
$$
\sum_{\textbf{n}=1}^{\infty}|\textbf{n}|^{-2}\sum_{\textbf{k} \preceq \textbf{n}}E\big[|X_{\textbf{k}}|I(|X_{\textbf{k}}|>|\textbf{n}|)\big] < \infty.$$\\

\begin{proof}
To prove the lemma, we take  $p=1$ in 
Lemma \ref{lemma2} and follow the proof up to  equation (\ref{L4.2E4}).  So, to prove our lemma, we have only to show that the series 
\begin{align}
\sum_{j=1}^{\infty}j^{-1}\int_{j}^{\infty}(\log^+ x)^{d-1}P(|X|>x)\,dx \label{L4.3E1}
\end{align}
is finite. \\

By comparing (\ref{L4.3E1}) with the integral, due to monotonicity, and by doing the change of order of integration, we obtain

\begin{align}
\bigintssss_{1}^{\infty}y^{-1}\bigintssss_{y}^{\infty}(\log^+ x)^{d-1}P(|X|>x)\,dx\, dy  
&= \bigintssss_{1}^{\infty}(\log^+ x)^{d-1} P(|X|>x)\bigintssss_{1}^{x}y^{-1} \,dy \, dx \notag \\
&\ll \bigintssss_{1}^{\infty}(\log^+ x)^{d}P(|X|>x) \, dx \notag  \\
&\ll E\big[|X|(\log^+|X|)^{d}\big] <\infty. \notag
\end{align}

Thus, series (\ref{L4.3E1}) is finite. Now, by following the proof of Lemma \ref{lemma2}, we obtain
$$\sum_{\textbf{n}=1}^{\infty}|\textbf{n}|^{-2}\sum_{\textbf{k} \preceq \textbf{n}}E\big[|X_{\textbf{k}}|I(|X_{\textbf{k}}|>|\textbf{n}|)\big] < \infty.$$
\end{proof}
\end{lemma}

\begin{lemma}\label{lemma4}
Assume $(X_{\textbf{n}})_{\textbf{n}}$ satisfies condition (MD). Let $1\leq p<2$, $1\leq r\leq p$.\\

If $E|X|^p < \infty$, 
then 
$$(i) 
\sum_{n=1}^{\infty}|N|^{(p-1)/r-1-1/d}\sum_{\textbf{k} \preceq N}E\big[|X_{\textbf{k}}|I(|X_{\textbf{k}}|>|N|^{1/r})\big] < \infty, \quad \text{for } p>1; 
$$
and,
$$(ii)
\sum_{n=1}^{\infty}|N|^{(p-2)/r-1-1/d}\sum_{\textbf{k} \preceq N}E\big[|X_{\textbf{k}}|^2I(|X_{\textbf{k}}|\leq|N|^{1/r})\big] < \infty, \quad \text{for } p\ge1.
$$
\begin{proof}
(i) Using (MD) and arguments similar to Lemma \ref{lemma2}, we notice that
\begin{align}
\sum_{n=1}^{\infty}|N|^{(p-1)/r-1-1/d}\sum_{\textbf{k} \preceq N}E\big[|X_{\textbf{k}}|I(|X_{\textbf{k}}|>|N|^{1/r})\big] 
&\ll\sum_{n=1}^{\infty}|N|^{p/r-1/d}P(|X|>|N|^{1/r}) \notag \\ 
&+\sum_{n=1}^{\infty}|N|^{(p-1)/r-1/d}\bigintssss_{|N|^{1/r}}^{\infty}P(|X|>x)\,dx. \label{L4.4E1}
\end{align}
Notice that $|N|=n^d$. Now, rewriting both the series on the right side of (\ref{L4.4E1}), we get 
\begin{align}
I= \sum_{n=1}^{\infty}n^{pd/r-1}  P(|X|>n^{d/r}) 
\leq E|X|^p, \label{L4.4E2}
\end{align}
and
\begin{align}
II &= \sum_{n=1}^{\infty}n^{(p-1)d/r-1}\bigintssss_{n^{d/r}}^{\infty}P(|X|>x)\,dx 
=\sum_{i=1}^{\infty}\bigintssss_{i^{d/r}}^{(i+1)^{d/r}}P(|X|>x)\,dx \sum_{n=1}^{i}n^{(p-1)d/r-1} \notag \\
&\ll \sum_{i=1}^{\infty}i^{(p-1)d/r}\bigintssss_{i^{d/r}}^{(i+1)^{d/r}} P(|X|>x)\,dx 
\leq \sum_{i=1}^{\infty}\bigintssss_{i^{d/r}}^{(i+1)^{d/r}} x^{p-1} P(|X|>x)\,dx  
\leq E|X|^p. \label{L4.4E3}
\end{align}

Combining (\ref{L4.4E1}), (\ref{L4.4E2}) and (\ref{L4.4E3}) gives
\begin{align}
\sum_{n=1}^{\infty}|N|^{(p-1)/r-1-1/d}\sum_{\textbf{k} \preceq N}E\big[|X_{\textbf{k}}|I(|X_{\textbf{k}}|>|N|^{1/r})\big] 
&\ll E|X|^p < \infty. \notag  
\end{align}

(ii) Condition (MD) gives
\begin{align}
\sum_{n=1}^{\infty}|N|^{(p-2)/r-1-1/d}\sum_{\textbf{k} \preceq N}E\big[|X_{\textbf{k}}|^2I(|X_{\textbf{k}}|\leq|N|^{1/r})\big]
&\ll \sum_{n=1}^{\infty}|N|^{(p-2)/r-1/d}\bigintssss_{0}^{|N|^{2/r}}P(X^2>x)\,dx.  \label{L4.4.2E1}  
\end{align}

Writing $\bigintss_{0}^{|N|^{2/r}}P(X^2>x)\,dx 
= \sum_{i=1}^{|N|}\bigintss_{(i-1)^{2/r}}^{i^{2/r}}P(X^2>x)\,dx$, and doing the change of order of summation 
\begin{align}
\sum_{n=1}^{\infty}|N|^{(p-2)/r-1/d}\bigintssss_{0}^{|N|^{2/r}}P(X^2>x)\,dx
&\ll 
\sum_{i=1}^{\infty}\bigintssss_{(i-1)^{2/r}}^{i^{2/r}}P(X^2>x)\,dx\sum_{n=i^{1/d}}^{\infty}n^{(p-2)d/r-1}  \notag \\
\ll  \sum_{i=1}^{\infty}i^{(p-2)/r}\bigintssss_{(i-1)^{2/r}}^{i^{2/r}}P(X^2>x)\,dx 
&\leq \sum_{i=1}^{\infty}\bigintssss_{(i-1)^{2/r}}^{i^{2/r}}x^{p/2-1}P(X^2>x)\,dx 
\ll E|X|^p. \label{L4.4.2E2}
\end{align}

Thus, (\ref{L4.4.2E1}) and (\ref{L4.4.2E2}) together imply that
\begin{align*}
\sum_{n=1}^{\infty}|N|^{(p-2)/r-1-1/d}\sum_{\textbf{k} \preceq N}E\big[|X_{\textbf{k}}|^2I(|X_{\textbf{k}}|\leq|N|^{1/r})\big]
\ll E|X|^p < \infty.  
\end{align*}
\end{proof}
\end{lemma}

\begin{lemma}\label{lemma5}
Assume $(X_{\textbf{n}})_{\textbf{n}}$ satisfies condition (MD)
and $E|X|log^+|X| < \infty$. 
Then 
$$
\sum_{n=1}^{\infty}|N|^{-1/d-1}\sum_{\textbf{k} \preceq N}E\big[|X_{\textbf{k}}|I(|X_{\textbf{k}}|>|N|)\big] < \infty.$$
\begin{proof}
The proof is similar to that of Lemma \ref{lemma3}, and therefore omitted.\\    
\end{proof}
\end{lemma}

\section{Proof of main results}\label{theoemproofs}

In this section, we now demonstrate the proof of our main results.

\subsection{Proof of \textit{Theorem \ref{Theorem1}}}
\begin{proof}
Let 
$$Y_\textbf{k}^{'}=e^{i\textbf{k}\cdot t}X_{\textbf{k}}I(|X_\textbf{k}|\leq |\textbf{n}|), \quad Y_\textbf{k}^{''}=e^{i\textbf{k}\cdot t}X_{\textbf{k}}I(|X_\textbf{k}|>|\textbf{n}|).$$

Note that
$S_{\textbf{n}}(t)=S_{\textbf{n}}^{'}(t)+S_{\textbf{n}}^{''}(t), $
where
$S_{\textbf{n}}^{'}(t)= \sum_{\textbf{k}\preceq \textbf{n}} Y_\textbf{k}^{'} \text{ and } S_{\textbf{n}}^{''}(t)= \sum_{\textbf{k}\preceq \textbf{n}} Y_\textbf{k}^{''}.$\\
        
Using Fubini's theorem, fact (\ref{EXlowertail}) and (MD) give
\begin{align}
\bigintssss_{\tilde{\Omega}}\frac{|S_{\textbf{n}}^{'}(t)|^2}{|\textbf{n}|^2}\,d\tilde{P}
&= \frac{1}{|\textbf{n}|^2}E\bigintssss_{\mathbb{T}^d}|S_{\textbf{n}}^{'}(t)|^2\,L(dt) 
=\frac{(2\pi)^d}{|\textbf{n}|^2}E\bigg(\sum_{\textbf{k}\preceq\textbf{n}}X_{\textbf{k}}^2I(|X_\textbf{k}|\leq |\textbf{n}|)\bigg) \notag\\
&\leq \frac{(2\pi)^d}{|\textbf{n}|^2}\sum_{\textbf{k}\preceq\textbf{n}}\int_{0}^{|\textbf{n}|^2}P(X_{\textbf{k}}^2>x) dx
\ll
\frac{1}{|\textbf{n}|}\int_{0}^{|\textbf{n}|^2}P(X^2>x) dx. \label{T1E1}
\end{align}

Now, writing $\int_{0}^{|\textbf{n}|^2}P(X^2>x) dx$ as $\sum_{i=1}^{|\textbf{n}|}\bigintssss_{(i-1)^2}^{i^2}P(X^2>x)\,dx$, note that the series
$\sum_{i=1}^{\infty} i^{-1}  \bigintssss_{(i-1)^2}^{i^2}P(X^2>x)\,dx
\leq  \bigintssss_{0}^{\infty}x^{-1/2}P(X^2>x)\,dx
\ll E|X| <\infty$. By Kronecker's lemma, this implies that the 
sequence on the right of (\ref{T1E1}) approaches to $0$ as $\max\{n_1,\dots,n_d\} \rightarrow \infty$. 

This implies that $\bigintssss_{\tilde{\Omega}}\frac{|S_{\textbf{n}}^{'}(t)|^2}{|\textbf{n}|^2}\,d\tilde{P}$ approaches to $0$ as $\max\{n_1,\dots,n_d\} \rightarrow \infty$. Therefore, we obtain
\begin{align}
\lim_{\max\{n_1,\dots,n_d\} \rightarrow \infty}\frac{S_{\textbf{n}}^{'}(t)}{|\textbf{n}|}=0 \quad \text{in the }L_1(\tilde{\Omega})\text{-norm}.\label{T1E2}
\end{align}
    
Next, by using Fubini's theorem, 
fact (\ref{EXuppertail}) and (MD), we obtain 
\begin{align*}
\bigintssss_{\tilde{\Omega}}\frac{|S_{\textbf{n}}^{''}(t)|}{|\textbf{n}|}\,d\tilde{P}
= \frac{1}{|\textbf{n}|}E\bigintssss_{\mathbb{T}^d}|S_{\textbf{n}}^{''}(t)|\,L(dt) 
\ll\frac{1}{|\textbf{n}|}E\bigg(\sum_{\textbf{k}\preceq\textbf{n}}|X_{\textbf{k}}|I(|X_\textbf{k}|> |\textbf{n}|)\bigg) 
\ll E|X|I(|X|> |\textbf{n}|).
\end{align*}

By the monotone convergence theorem, $E|X|I(|X|> |\textbf{n}|) \rightarrow 0$ as $\max\{n_1,\dots,n_d\} \rightarrow \infty$, and therefore
\begin{align}
\lim_{\max\{n_1,\dots,n_d\} \rightarrow \infty}\frac{S_{\textbf{n}}^{''}(t)}{|\textbf{n}|}=0 \quad \text{in the }L_1(\tilde{\Omega})\text{-norm}.\label{T1E3}
\end{align}

Combining (\ref{T1E2}) and (\ref{T1E3}), we obtain
$$
\lim_{\max\{n_1,\dots,n_d\} \rightarrow \infty}\frac{S_{\textbf{n}}(t)}{|\textbf{n}|}=0 \quad \text{in the }L_1(\tilde{\Omega})\text{-norm}.
$$

\end{proof}

\subsection{Proof of Theorem \ref{Theorem2}}

\begin{proof}
For $t \in \mathbb{T}^d$, define the truncated random fields

$$X_\textbf{k}^{'}=e^{i\textbf{k}\cdot \textbf{t}}X_{\textbf{k}}I(|X_\textbf{k}|\leq |\textbf{n}|^{1/r}), \quad X_\textbf{k}^{''}=e^{i\textbf{k}\cdot \textbf{t}}X_{\textbf{k}}I(|X_\textbf{k}|>|\textbf{n}|^{1/r}),$$

and their partial sums
$$
S_{\textbf{n}}^{'}(\textbf{t})= \sum_{\textbf{k}\preceq \textbf{n}} X_\textbf{k}^{'} \text{ and } S_{\textbf{n}}^{''}(\textbf{t})= \sum_{\textbf{k}\preceq \textbf{n}} X_\textbf{k}^{''}.$$

Clearly,
$$e^{i\textbf{k}\cdot \textbf{t}}X_\textbf{k}=X_\textbf{k}^{'}+X_\textbf{k}^{''} \text{ and } 
S_{\textbf{n}}(\textbf{t})=S_{\textbf{n}}^{'}(\textbf{t})+S_{\textbf{n}}^{''}(\textbf{t}). $$
\\

By the Markov's Inequality,

\begin{equation}
 P\big(|S_{\textbf{n}}^{'}(t)|>\frac{\epsilon |\textbf{n}|^{1/r}}{2}\big) \ll |\textbf{n}|^{-2/r}E[|S_{\textbf{n}}^{'}(t)|^2].\label{T2E1}
\end{equation}

Then by using Fubini's theorem,
\begin{align}
 \tilde{P}\big(|S_{\textbf{n}}^{'}(t)|>\frac{\epsilon |\textbf{n}|^{1/r}}{2}\big)
& =\int_{\mathbb{T}^d} P\big(|S_{\textbf{n}}^{'}(t)|>\frac{\epsilon |\textbf{n}|^{1/r}}{2}\big)\,L(dt) 
\ll |\textbf{n}|^{-2/r} \int_{\mathbb{T}^d} E[|S_{\textbf{n}}^{'}(t)|^2]\,L(dt) \notag \\
&\ll |\textbf{n}|^{-2/r} E\bigg[\sum_{\textbf{k}\preceq \textbf{n}} |X_\textbf{k}^{'}|^2\bigg] 
=|\textbf{n}|^{-2/r}\sum_{\textbf{k}\preceq \textbf{n}} E\big[|X_\textbf{k}|^2I(|X_\textbf{k}|\leq |\textbf{n}|^{1/r})\big]. \label{T2E2}
\end{align}

Now, Lemma \ref{lemma2} gives
\begin{equation}
\sum_{\textbf{n}=1}^{\infty}|\textbf{n}|^{p/r-2} \tilde{P}\big(|S_{\textbf{n}}^{'}(t)|>\frac{\epsilon |\textbf{n}|^{1/r}}{2}\big)  
\ll \sum_{\textbf{n}=1}^{\infty}|\textbf{n}|^{(p-2)/r-2}\sum_{\textbf{k}\preceq \textbf{n}} E\big[|X_\textbf{k}|^2I(|X_\textbf{k}|\leq |\textbf{n}|^{1/r})\big] < \infty. \label{T2E3}
\end{equation}
\\

By the Markov's Inequality 

\begin{equation}
P\big(|S_{\textbf{n}}^{''}(t)|>\frac{\epsilon |\textbf{n}|^{1/r}}{2}\big) 
 \ll |\textbf{n}|^{-1/r}E[|S_{\textbf{n}}^{''}(t)]\\
=|\textbf{n}|^{-1/r}\sum_{\textbf{k}\preceq \textbf{n}}E\big[|X_\textbf{k}|I(|X_\textbf{k}|>|\textbf{n}|^{1/r})\big]. \label{T2E4}
\end{equation}

Now, by using Fubini's theorem along with (\ref{T2E4}) and Lemma \ref{lemma2}, we get
\begin{align}
\sum_{\textbf{n}=1}^{\infty}|\textbf{n}|^{p/r-2} \tilde{P}\big(|S_{\textbf{n}}^{''}(t)|>\frac{\epsilon |\textbf{n}|^{1/r}}{2}\big)
& =\sum_{\textbf{n}=1}^{\infty}|\textbf{n}|^{p/r-2}\int_{\mathbb{T}^d} P\big(|S_{\textbf{n}}^{''}(t)|>\frac{\epsilon |\textbf{n}|^{1/r}}{2}\big)\,L(dt) \notag \\
&\ll\sum_{\textbf{n}=1}^{\infty}|\textbf{n}|^{(p-1)/r-2}\int_{\mathbb{T}^d}\sum_{\textbf{k}\preceq \textbf{n}}E\big[|X_\textbf{k}|I(|X_\textbf{k}|>|\textbf{n}|^{1/r})\big] \,L(dt) \notag \\
&=(2\pi)^d \sum_{\textbf{n}=1}^{\infty}|\textbf{n}|^{(p-1)/r-2}\sum_{\textbf{k}\preceq \textbf{n}}E\big[|X_\textbf{k}|I(|X_\textbf{k}|>|\textbf{n}|^{1/r})\big] < \infty. \label{T2E5}
\end{align}
\\

Combining  (\ref{T2E3}) and (\ref{T2E5}) gives

\begin{align}
\sum_{\textbf{n}=1}^{\infty}|\textbf{n}|^{p/r-2} \tilde{P}\big(|S_{\textbf{n}}(t)|>\epsilon|\textbf{n}|^{1/r}\big)
\leq 
\sum_{\textbf{n}=1}^{\infty}|\textbf{n}|^{p/r-2} \tilde{P}\big(|S_{\textbf{n}}^{'}(t)|>\frac{\epsilon |\textbf{n}|^{1/r}}{2}\big)
+
\sum_{\textbf{n}=1}^{\infty}|\textbf{n}|^{p/r-2} \tilde{P}\big(|S_{\textbf{n}}^{''}(t)|>\frac{\epsilon |\textbf{n}|^{1/r}}{2}\big)
< \infty. \label{T2E6}
\end{align}
\\

Finally, by using Fubini's theorem together with (\ref{T2E6}), we obtain 
$$
\int_{\mathbb{T}^d} \sum_{\textbf{n}=1}^{\infty}|\textbf{n}|^{p/r-2}P\big(|S_{\textbf{n}}(t)|>\epsilon |\textbf{n}|^{1/r}\big)\,L(dt)
=
\sum_{\textbf{n}=1}^{\infty}|\textbf{n}|^{p/r-2} \tilde{P}\big(|S_{\textbf{n}}(t)|>\epsilon|\textbf{n}|^{1/r}\big)
< \infty.
$$
\\

This shows that for almost all $t \in \mathbb{T}^d$
$$
\sum_{\textbf{n}=1}^{\infty}|\textbf{n}|^{p/r-2} P\big(|S_{\textbf{n}}(t)|>\epsilon|\textbf{n}|^{1/r}\big)< \infty.
$$
\end{proof}

\subsection{Proof of Theorem \ref{Theorem3}}
\begin{proof}
By taking $p=1$ and using arguments similar to the proof of Theorem \ref{Theorem2}, we obtain
\begin{align}
 \tilde{P}\big(|S_{\textbf{n}}^{'}(t)|>\frac{\epsilon |\textbf{n}|}{2}\big)   
 \ll |\textbf{n}|^{-2}\sum_{\textbf{k}\preceq \textbf{n}} E\big[|X_\textbf{k}|^2I(|X_\textbf{k}|\leq |\textbf{n}|)\big], \label{T3E1}  
\end{align}
and
\begin{align}
\tilde{P}\big(|S_{\textbf{n}}^{''}(t)|>\frac{\epsilon |\textbf{n}|}{2}\big)
\ll 
|\textbf{n}|^{-1}\sum_{\textbf{k}\preceq \textbf{n}}E\big[|X_\textbf{k}|I(|X_\textbf{k}|>|\textbf{n}|)\big].\label{T3E2}
\end{align}
\\

Then, using Lemma \ref{lemma2} and Lemma \ref{lemma3}, we obtain
\begin{align}
\sum_{\textbf{n}=1}^{\infty}|\textbf{n}|^{-1} \tilde{P}\big(|S_{\textbf{n}}^{'}(t)|>\frac{\epsilon |\textbf{n}|}{2}\big)  
\ll \sum_{\textbf{n}=1}^{\infty}|\textbf{n}|^{-3}\sum_{\textbf{k}\preceq \textbf{n}} E\big[|X_\textbf{k}|^2I(|X_\textbf{k}|\leq |\textbf{n}|)\big] < \infty, \label{T3E3}    
\end{align}
and
\begin{align}
\sum_{\textbf{n}=1}^{\infty}|\textbf{n}|^{-1} \tilde{P}\big(|S_{\textbf{n}}^{''}(t)|>\frac{\epsilon |\textbf{n}|}{2}\big)  
\ll \sum_{\textbf{n}=1}^{\infty}|\textbf{n}|^{-2}\sum_{\textbf{k}\preceq \textbf{n}} E\big[|X_\textbf{k}|I(|X_\textbf{k}|>|\textbf{n}|)\big] < \infty. \label{T3E4}
\end{align}
\\

Combining (\ref{T3E3}) and (\ref{T3E4}), we get
\begin{align}
\sum_{\textbf{n}=1}^{\infty}|\textbf{n}|^{-1} \tilde{P}\big(|S_{\textbf{n}}(t)|>\frac{\epsilon |\textbf{n}|}{2}\big)  
< \infty. \label{T3E5}
\end{align}
\\

Now, (\ref{T3E5}) along with using Fubini's theorem gives, for almost all $t \in \mathbb{T}^d$  
$$
\sum_{\textbf{n}=1}^{\infty}|\textbf{n}|^{-1} P\big(|S_{\textbf{n}}(t)|>\frac{\epsilon |\textbf{n}|}{2}\big)  
< \infty.$$

\end{proof}

\subsection{Proof of Theorem \ref{Theorem4}}
\begin{proof}
For $t \in \mathbb{T}^d$, define the truncated random fields as defined in Theorem \ref{Theorem2} over a cubic domain $\textbf{k} \preceq \textbf{n}$, where we take $\textbf{n}=N=(n,\dots,n) \in \mathbb{Z}_+^d$.\\

The proof of our theorem is similar to that of Theorem \ref{Theorem2}.\\

By using Fubini's theorem together with the Markov's Inequality and the maximal inequality in Theorem 4.4 of \cite{weisz2012summability}, we get
\begin{align}
 \tilde{P}\bigg(\underset{1 \leq m \leq n}{\max}
|S^{'}_{M}(t)|>\frac{\epsilon|N|^{1/r}}{2}\bigg)
&\ll |N|^{-2/r} \int_{\mathbb{T}^d} E\bigg[\bigg(\underset{1 \leq m \leq n}{\max}|S^{'}_{M}(t)|\bigg)^2 \bigg]\,L(dt) \notag \\
&\ll |N|^{-2/r} E\bigg[\sum_{\textbf{k}\preceq N} |X_\textbf{k}^{'}|^2\bigg] \notag \\
&=|N|^{-2/r}\sum_{\textbf{k}\preceq N} E\big[|X_\textbf{k}|^2I(|X_\textbf{k}|\leq |N|^{1/r})\big]. \label{T4E1}
\end{align}

Now, Lemma \ref{lemma4} gives
\begin{equation}
\sum_{n=1}^{\infty}|N|^{p/r-1-1/d} \tilde{P}\bigg(\underset{1 \leq m \leq n}{\max}
|S^{'}_{M}(t)|>\frac{\epsilon|N|^{1/r}}{2}\bigg) 
\ll \sum_{n=1}^{\infty}|N|^{(p-2)/r-1-1/d}\sum_{\textbf{k}\preceq N} E\big[|X_\textbf{k}|^2I(|X_\textbf{k}|\leq |N|^{1/r})\big] < \infty. \label{T4E2}\\
\end{equation}
\\

Also, by using Fubini's theorem along with the Markov's Inequality and Lemma \ref{lemma4}, we get

\begin{align}
\sum_{n=1}^{\infty}|N|^{p/r-1-1/d} \tilde{P}\bigg(\underset{1 \leq m \leq n}{\max}
|S^{''}_{M}(t)|>\frac{\epsilon|N|^{1/r}}{2}\bigg)
&\ll \sum_{n=1}^{\infty}|N|^{(p-1)/r-1-1/d} \int_{\mathbb{T}^d} E\big[\underset{1 \leq m \leq n}{\max}|S^{''}_{M}(t)| \big]\,L(dt) \notag \\
&\ll \sum_{n=1}^{\infty}|N|^{(p-1)/r-1-1/d} \int_{\mathbb{T}^d}E\bigg[\sum_{\textbf{k}\preceq N} |X_\textbf{k}^{''}|\bigg] \,L(dt) \notag \\
&=(2\pi)^d \sum_{n=1}^{\infty}|N|^{(p-1)/r-1-1/d} \sum_{\textbf{k}\preceq N}E\big[|X_\textbf{k}|I(|X_\textbf{k}|>|N|^{1/r})\big] 
< \infty. \label{T4E3}
\end{align}

Combining  (\ref{T4E2}) and (\ref{T4E3}) gives

\begin{align}
\sum_{n=1}^{\infty}|N|^{p/r-1-1/d} \tilde{P}\bigg(\underset{1 \leq m \leq n}{\max}
|S_{M}(t)|>\epsilon|N|^{1/r}\bigg)
< \infty. \label{T4E5}
\end{align}
\\

Finally, by using Fubini's theorem together with (\ref{T4E5}), we obtain 
$$
\bigintsss_{\mathbb{T}^d} \sum_{n=1}^{\infty}|N|^{p/r-1-1/d} P\bigg(\underset{1 \leq m \leq n}{\max}
|S_{M}(t)|>\epsilon|N|^{1/r}\bigg)\, L(dt)
=
\sum_{n=1}^{\infty}|N|^{p/r-1-1/d} \tilde{P}\bigg(\underset{1 \leq m \leq n}{\max}
|S_{M}(t)|>\epsilon|N|^{1/r}\bigg)
< \infty.
$$
\\

This shows that for almost all $t \in \mathbb{T}^d$
\begin{align}
\sum_{n=1}^{\infty}|N|^{p/r-1-1/d}P\bigg(\underset{1 \leq m \leq n}{\max}
|S_{M}(t)|>\epsilon|N|^{1/r}\bigg)
< \infty. \notag   
\end{align}
\end{proof}

\subsection{Proof of Theorem \ref{Theorem5}}
\begin{proof}
The proof to this thoerem follows from that of Theorem \ref{Theorem4} except for the fact that we use Lemma \ref{lemma5} instead of Lemma \ref{lemma4} in (\ref{T4E3}).\\  
\end{proof}

\subsection{Proof of Remark \ref{remark1}}
Letting $r=p$ in Theorem \ref{Theorem4}, we obtain
$\sum_{n=1}^{\infty}|N|^{-1/d} P\bigg(\underset{1 \leq m \leq n}{\max}
|S_{M}(t)|>\epsilon|N|^{1/p}\bigg)
< \infty.$

Now, this implies that

\begin{align*}
\infty > &\sum_{n=1}^{\infty}|N|^{-1/d} P\bigg(\underset{1 \leq m \leq n}{\max}
|S_{M}(t)|>\epsilon|N|^{1/p}\bigg)   \\
&= \sum_{k=1}^{\infty}\sum_{n=2^{k-1}}^{2^k}n^{-1} P\bigg(\underset{1 \leq m \leq n}{\max}
|S_{M}(t)|>\epsilon n^{d/p}\bigg) \\
& \geq \frac{1}{2}\sum_{k=1}^{\infty}P\bigg(\underset{1\leq m\leq 2^{k-1}}{\max}
|S_{M}(t)|>\epsilon 2^{kd/p}\bigg). \\
\end{align*}

By using the Borel-Cantelli Lemma, this implies that for almost all $t \in \mathbb{T}^d$ 
$$2^{-kd/p}\underset{1\leq m\leq 2^{k}}{\max}
|S_{M}(t)| \xrightarrow{\text{P-a.s.}} 0 \quad \text{as } k \rightarrow \infty,$$

which implies that
$$\frac{S_{N}(t)}{|N|^{1/p}} \xrightarrow{\text{P-a.s.}} 0 \quad \text{as } n \rightarrow \infty.$$\\

\section{Appendix}

\begin{result}
   By \cite{smythe1974sums}, we have
   \begin{equation}\label{eq:1}
       M(x)=O(x(\log^+x)^{d-1}) \quad \text{and} \quad d(x)=o(x^{\delta}), \quad  \forall \delta >0 \text{ as } x\rightarrow \infty.
   \end{equation}\\
\end{result}

\begin{result}
Following are the known facts: For $X \geq 0$ a.s. and $\alpha > 0$,
\begin{equation}\label{EXuppertail}
E\big[XI(X>\alpha)\big]=\alpha P(X>\alpha)+\int_{\alpha}^{\infty}P(X>x)dx,
\end{equation}
and
\begin{equation}\label{EXlowertail}
E\big[XI(X\leq\alpha)\big]=-\alpha P(X>\alpha)+\int_{0}^{\alpha}P(X>x)dx.    
\end{equation}\\
\end{result}

\begin{result}
We have the following elementary inequalities:
\begin{equation}
   \sum_{j=i}^{\infty}j^{-\alpha} \ll i^{1-\alpha} \quad \text{for }\alpha>1. 
\end{equation} 

\begin{equation}
   \sum_{j=1}^{i}j^{\alpha} \ll i^{1+\alpha} \quad \text{for }\alpha>-1. 
\end{equation}\\
\end{result}

\begin{result}
Let $f$ be an increasing function on $\mathbb{R}$. Then, for any $\alpha,\beta \ne 0$ 
\begin{equation}\label{EXf(X)}
E[|X|^{\alpha}f(|X|)] \geq \frac{\alpha}{\beta}\int_{0}^{\infty}x^{\frac{\alpha}{\beta}-1}f(x^{1/{\beta}})P(|X|^{\beta}>x) \,dx.   
\end{equation}
\begin{proof}
Note that $|X|^{\alpha}=\frac{\alpha}{\beta}\bigintsss_{0}^{\infty} x^{\frac{\alpha}{\beta}-1}I(|X|^{\beta}>x) \, dx$. 
Then,
\begin{align*}
E[|X|^{\alpha}f(|X|)]
&= \int |X|^{\alpha}f(|X|) \,dP
=\frac{\alpha}{\beta}\int_{\Omega} \bigg(\int_{0}^{\infty} x^{\frac{\alpha}{\beta}-1}I(|X(w)|^{\beta}>x)\, dx\bigg) f(|X(w)|) \,P(dw)\\
&=\frac{\alpha}{\beta}\int_{0}^{\infty} x^{\frac{\alpha}{\beta}-1}\bigg(\int_{|X(w)|^{\beta}>x}f(|X(w)|) \,P(dw)\bigg) \, dx\\
&\geq\frac{\alpha}{\beta}\int_{0}^{\infty} x^{\frac{\alpha}{\beta}-1}f(x^{\frac{1}{\beta}})\bigg(\int_{|X(w)|^{\beta}>x} \,P(dw)\bigg) \, dx\\
&=\frac{\alpha}{\beta}\int_{0}^{\infty} x^{\frac{\alpha}{\beta}-1}f(x^{\frac{1}{\beta}})P(|X|^{\beta}>x) \, dx.\\
\end{align*}
\end{proof}
\end{result}

\section*{Acknowledgments}
This research was supported by the NSF grant: DMS-2054598 and DMS-2336047.

The author would like to thank Magda Peligrad for providing her insightful guidance in the completion of this work.

\bibliographystyle{plain}  
\bibliography{ref}

\end{document}